\documentclass[11pt]{article}
\usepackage{latexsym,amsmath,stackrel,color,amsthm,amssymb,epsfig,graphicx,mathrsfs}
\usepackage{graphicx,comment}
\usepackage{amssymb}
\usepackage[left=1in,top=1in,right=1in,bottom=1in]{geometry}
\usepackage[linktocpage=true]{hyperref}
\usepackage{setspace}
\usepackage{amssymb, amsmath, amsthm, graphicx,mathrsfs}
\usepackage{caption}

\usepackage{tikz}
\makeatletter
\def\thm@space@setup{%
  \thm@preskip=\parskip \thm@postskip=0pt
}
\makeatother
\def\qed{\hfill\ifhmode\unskip\nobreak\fi\quad\ifmmode\Box\else\hfill$\Box$\fi}
\def\ite#1{\hfill\break${}$\hbox to 50pt {\quad(#1)\hfill}}

\newtheorem{thm}{Theorem}
\newtheorem{cor}[thm]{Corollary}

\def\ex{{\rm{ex}}}
\def\ord{{\rm{ord}}}

\begin{document}

\pagestyle{myheadings}
\markright{{\small{\sc  Mubayi and Verstra\"ete:   Ramsey numbers
}}}

\title{
\vspace{-0.9in}A note on pseudorandom Ramsey graphs}

\author{
Dhruv Mubayi\thanks{Research partially supported by NSF award DMS-1763317.} \and Jacques Verstra\"ete\thanks{Research partially supported by NSF award DMS-1800832. \newline
2010 MSC: 05D10, 05C55, 05B25}
}

\maketitle

\begin{abstract} For fixed $s \ge 3$, we prove that if optimal $K_s$-free pseudorandom graphs exist, then the Ramsey number
	$r(s,t) = t^{s-1+o(1)}$ as $t \rightarrow \infty$. Our method also improves the best lower bounds for $r(C_{\ell},t)$ obtained by Bohman and Keevash from
the random $C_{\ell}$-free process by polylogarithmic factors for all odd $\ell \geq 5$ and $\ell \in \{6,10\}$. For $\ell = 4$ it matches their
lower bound from the $C_4$-free process.

We also prove, via a different approach,  that $r(C_5, t)> (1+o(1))t^{11/8}$ and $r(C_7, t)> (1+o(1))t^{11/9}$. These improve the exponent of $t$ in the previous best results and appear to be the first examples of graphs $F$ with cycles for which such an improvement of the exponent for $r(F, t)$ is shown over the bounds given by the random $F$-free process and random graphs.
\end{abstract}

\section{Introduction}

The {\em Ramsey number} $r(F,t)$ is the minimum $N$ such that every $F$-free graph on $N$ vertices has an independent set of size $t$. When $F = K_s$ we simply write $r(s,t)$ instead of $r(F,t)$. Improving on earlier results of Spencer~\cite{S} and the classical Erd\H os-Szekeres~\cite{ES} theorem on Ramsey numbers, Ajtai, Koml\'{o}s and Szemer\'{e}di~\cite{AKS}
proved the following upper bound on $r(s,t)$, and Bohman and Keevash~\cite{BK2} proved the lower bound by considering the {\em random $K_s$-free process}: consequently
for $s \ge 3$, there exist constants $c_1(s),c_2(s) > 0$ such that
\begin{equation} \label{rsn}
c_1(s) \frac{t^{\frac{s+1}{2}}}{(\log t)^{\frac{s + 1}{2} - \frac{1}{s-2}}} \leq r(s,t) \leq c_2(s) \frac{t^{s-1}}{(\log t)^{s - 2}}.
\end{equation}
For $s = 3$, the lower bound was proved in a celebrated paper of Kim~\cite{Kim} and the upper bound was proved by Shearer~\cite{Sh} with $c_2(3) = 1 + o(1)$.
In particular, recent results of Bohman and Keevash~\cite{BK} and Fiz Pontiveros, Griffiths and Morris~\cite{PGM} together with the bound of Shearer show
\begin{equation}
(\tfrac{1}{4} - o(1)) \cdot \frac{t^2}{\log t} \leq r(3,t) \leq (1 + o(1)) \cdot \frac{t^2}{\log t}
\end{equation}
as $t \rightarrow \infty$. There have been no improvements in the exponents in (\ref{rsn}) for any $s \ge 4$ for many decades.
In this note, we show that if certain density-optimal $K_s$-free pseudorandom graphs exist, then $r(s,t) = t^{s  - 1 + o(1)}$. This approach suggests that pseudorandom graphs may be the central tool required to determine classical graph Ramsey numbers. 

\medskip

An {\em $(n,d,\lambda)$ graph} is an $n$-vertex $d$-regular graph such that the absolute value of every eigenvalue of its adjacency matrix, besides the largest one, is at most $\lambda$. Constructions of $(n,d,\lambda)$-graphs arise from a number of sources, including Cayley graphs, projective geometry and strongly regular graphs -- we refer
the reader to Krivelevich and Sudakov~\cite{KS} for a survey of $(n,d,\lambda)$-graphs.  Sudakov, Szabo and Vu~\cite{SSV} show that a $K_s$-free $(n,d,\lambda)$-graph
satisfies
\begin{equation}\label{lambdabound}
\lambda = \Omega(d^{s-1}/n^{s-2})
\end{equation}
as $n \rightarrow \infty$. For $s = 3$, if $G$ is any triangle-free $(n,d,\lambda)$-graph with adjacency matrix $A$,  then
\begin{equation}\label{hoffman}
0 = \hbox{tr}(A^3) \ge  d^3 -\lambda^3(n-1).
\end{equation}
If $\lambda = O(\sqrt d)$, then this gives $d = O(n^{2/3})$ matching (\ref{lambdabound}).
 Alon~\cite{A} constructed a triangle-free pseudorandom graph attaining this bound, and Conlon~\cite{Conlon} more recently analyzed a randomized construction with the same average degree.
A similar argument to (\ref{hoffman}) shows that a $K_s$-free $(n,d,\lambda)$-graph with $\lambda = O(\sqrt{d})$ has $d = O(n^{1 - \frac{1}{2s-3}})$.
The Alon-Boppana Bound~\cite{Nilli1,Nilli2} shows that $\lambda = \Omega(\sqrt{d})$ for every $(n,d,\lambda)$-graph provided $d/n$ is bounded away from 1. Sudakov, Szabo and Vu~\cite{SSV} raised the question of the existence of optimal pseudorandom $K_s$-free graphs for $s \geq 4$, namely $(n,d,\lambda)$-graphs achieving the bound in (\ref{lambdabound}) with $\lambda = O(\sqrt{d})$ and $d = \Omega(n^{1 - \frac{1}{2s-3}})$.
We show that a positive answer to this question gives the exponent of the Ramsey numbers $r(s,t)$ via a short proof of the following general theorem, based on ideas of Alon and R\"{o}dl~\cite{AR}:

\begin{thm} \label{main}
Let $F$ be a graph, $n,d,\lambda$ be positive integers with $d\ge 1$ and $\lambda>1/2$ and let $t = \lceil 2 n\log^2\!n/d\rceil$. If there exists an $F$-free $(n,d,\lambda)$-graph, then
\begin{equation}\label{first}
r(F,t) >  \frac{n}{20 \lambda} \log^2\!n .
\end{equation}
\end{thm}

Theorem~\ref{main}  provides good bounds whenever we have an $F$-free $(n,d,\lambda)$-graph with many edges and good pseudorandom properties (meaning that $d$ is large and $\lambda$ is small). For example, we immediately obtain the following consequence.

\begin{cor} \label{cor1}
If  $K_s$-free  $(n,d,\lambda)$-graphs  exist with $d = \Omega(n^{1 - \frac{1}{2s-3}})$ and $\lambda = O(\sqrt{d})$, then as $t \rightarrow \infty$,
\begin{equation}\label{second}
r(s,t) = \Omega\Bigl(\frac{t^{s-1}}{\log^{2s-4} t}\Bigr).
\end{equation}
\end{cor}

Corollary~\ref{cor1}
follows from (\ref{first}) using $F = K_s$. To see this, from $t = \lceil 2n\log^2\!n/d \rceil$ we obtain $d = \Theta(t^{2s-4}/(\log t)^{2(2s-4)})$.
Inserting this in (\ref{first}) with $\lambda = O(\sqrt{d})$ gives (\ref{second}).

 Alon and Krivelevich~\cite{AKr}  gave a construction of $K_s$-free $(n,d,\lambda)$-graphs with $d = \Omega(n^{1 - 1/(s - 2)})$ and $\lambda = O(\sqrt{d})$ for all $s \geq 3$, and this was slightly improved by Bishnoi, Ihringer and Pepe~\cite{BIP} to obtain $d = \Omega(n^{1 - 1/(s - 1)})$. This is the current record
for the degree of a $K_s$-free $(n,d,\lambda)$-graph with $\lambda = O(\sqrt{d})$. The problem of obtaining optimal $K_s$-free pseudorandom constructions in the sense (\ref{lambdabound}) with $\lambda = O(\sqrt{d})$ for $s \geq 4$ seems difficult and is considered to be a central open problem in pseudorandom graph theory. The problem of determining the growth rate of $r(s,t)$ is classical and much older, and it wasn't completely clear whether the upper or lower bound in (\ref{rsn}) was closer to the truth. Based on Theorem~\ref{main}, it seems reasonable to conjecture that if $s \geq 4$ is fixed, then $r(s,t) = t^{s - 1 + o(1)}$ as $t \rightarrow \infty$. 

\medskip

We next consider cycle-complete Ramsey numbers. The cycle complete Ramsey numbers $r(C_{\ell},t)$ appear to be very difficult to determine -- the best upper bounds are provided by Sudakov~\cite{Sudakov} for odd cycles and Caro, Li, Rousseau and Zhang~\cite{CaroLi} for even cycles.
The best lower bound for fixed $\ell\ge 4$ is
 \begin{equation} \label{eqbk}
 r(C_{\ell}, t) = \Omega\left(\frac{t^{(\ell-1)/(\ell-2)}}{\log t}\right)
 \end{equation}
 due to Bohman and Keevash~\cite{BK2} by analyzing the $C_{\ell}$-free process.
A generalization of the optimal triangle-free $(n,d,\lambda)$-graphs constructed by Alon~\cite{A} to optimal pseudorandom $C_{\ell}$-free graphs for odd $\ell \geq 5$
was given by Alon and Kahale~\cite{AKa}, and gives an $(n,d,\lambda)$-graph with $d = \Theta(n^{2/\ell})$ and $\lambda = O(\sqrt{d})$. Using this construction,
Theorem \ref{main} gives the following on odd-cycle complete Ramsey numbers, which gives a polylogarithmic improvement over (\ref{eqbk}):

\medskip

\begin{cor} \label{cor1a}
Let $\ell \geq 3$ be an odd integer. Then as $t \rightarrow \infty$,
\begin{equation}\label{1a}
r(C_{\ell},t) = \Omega\Bigl(\frac{t^{(\ell - 1)/(\ell - 2)}}{\log^{2/(\ell - 2)}t}\Bigr).
\end{equation}
\end{cor}

\medskip

Note when $\ell = 3$, this matches the lower bound of Spencer~\cite{S} from the local lemma. Applying Theorem \ref{main} when $F$ is bipartite can give lower bounds on $r(F,t)$ that are better than those obtained from the $F$-free process.
We can see this when $F= C_{\ell}$ and $\ell \in \{6,10\}$.
 To apply Theorem\ref{main} when $F = C_4$, we may consider polarity graphs of projective planes to be $(n,d,\lambda)$-graphs with $n = q^2 + q + 1$, $d = q + 1$ and $\lambda = \sqrt{q}$ (see~\cite{MW} for a  detailed study of independent sets in such graphs).   Theorem \ref{main} then gives $r(C_4,t) = \Omega(t^{3/2}/\log t)$ which matches (\ref{eqbk}). It is a wide open conjecture of Erd\H{o}s that $r(C_4,t) \leq t^{2 - \epsilon}$ for some $\epsilon > 0$.

  For $\ell\in \{6,10\}$ Theorem~\ref{main} provides results that exceed the previous best known bounds of (\ref{eqbk}) from the random $C_{\ell}$-free process.

\begin{cor} \label{cor2}
As $t \rightarrow \infty$, 	
$$
r(C_6, t) = \Omega\left(\frac{t^{5/4}}{\log^{1/2} t}\right) \qquad \hbox{ and } \qquad r(C_{10}, t) = \Omega\left(\frac{t^{9/8}}{\log^{1/4} t}\right).$$
	\end{cor}
These results are obtained by considering polarity graphs of {\em generalized quadrangles} and {\em generalized hexagons}. For certain prime powers $q$,
generalized quadrangles are $(n,d,\lambda)$-graphs with $n = q^3 + q^2 + q + 1$, $d = q + 1$ and $\lambda = \sqrt{2q}$, and generalized hexagons
are $(n,d,\lambda)$-graphs with $n = q^5 + q^4 + q^3 + q^2 + q + 1$, $d = q + 1$ and $\lambda = \sqrt{3q}$. For the existence of such graphs, we refer
the reader to Brouwer, Cohen and Neumaier~\cite{BCN} and Lazebnik, Ustimenko and Woldar~\cite{LUW}.
We can then apply Theorem~\ref{main} and obtain the desired result by (\ref{first}).



Our next result uses a completely different construction than that in Theorem~\ref{main}
for $F = C_5$ and $F = C_7$. For these two cases, we are able to improve the exponents in the lower bounds given by Corollary \ref{cor1a}. Our approach here is to use a random block construction. This idea was used in~\cite{DRR} and \cite{KMV} and also recently in~\cite{Conlon} to construct triangle-free pseudorandom graphs.

\begin{thm}\label{main2}
As $t \rightarrow \infty$,
\begin{eqnarray*}
r(C_5,t) &\geq& (1 + o(1))t^{11/8}\\
r(C_7,t) &\geq& (1 + o(1))t^{11/9}.
\end{eqnarray*}
\end{thm}

This appears to be the first instance of a graph $F$ containing cycles for which random graphs do not supply the right exponent for $r(F,t)$.

\section{Proof of Theorem \ref{main}}

The proof of Theorem \ref{main} uses the following property of independent sets in $(n,d,\lambda)$-graphs, due to Alon and R\"{o}dl~\cite{AR} (we give a slightly
stronger statement below):

\begin{thm} {\rm (Alon-R\"odl~\cite{AR}) } \label{alonrodl}
Let $G$ be an $(n,d,\lambda)$-graph with $d \geq 1$, $\lambda>1/2$, and let $t \geq  2n\log^2\!n/d$ be an integer. Then the number of independent sets of size $t$ in $G$ is at most $(\frac{2e^2 \lambda}{\log^2\!n})^t$.
\end{thm}

\proof Alon and R\"odl (Theorem 2.1 in~\cite{AR}) proved that  the number $Z$ of independent sets of size $t$ in $G$ is at most
$$\frac{1}{t!} {t \choose \ell} n^{\ell} \left(\frac{2\lambda n}{d}\right)^{t-\ell}$$
where $\ell = t/\log n$. Using ${t \choose \ell} \leq 2^t$ and $t! \geq (t/e)^t$,
$$Z \leq \Bigl(\frac{2e}{t}\Bigr)^t \cdot n^{\ell} \cdot \Bigl(\frac{2\lambda n}{d}\Bigr)^{t} \cdot \Bigl(\frac{d}{2\lambda n}\Bigr)^{\ell}
= \Bigl(\frac{4e\lambda n}{dt}\Bigr)^t \cdot \Bigl(\frac{d}{2\lambda}\Bigr)^{t/\log n}.$$
Since  $\lambda>1/2$, we obtain $d/2\lambda \leq d \leq n$ and therefore  $(d/2\lambda)^{t/\log n} \leq n^{t/\log n} \leq e^t$. Using $t \geq 2n(\log^2 n)/d$,
$$Z \leq \Bigl(\frac{4e^2 \lambda n}{dt}\Bigr)^t \leq \Bigl(\frac{2e^2 \lambda}{\log^2\!n}\Bigr)^t$$
and the proof is complete. \qed

\medskip

{\bf Proof of Theorem \ref{main}.} Let $G$ be an $F$-free $(n,d,\lambda)$-graph and let $U$ be a random set of vertices of $G$ where each vertex is chosen
 independently with probability $p =\log^2 n/2e^2\lambda$. Let $Z$ be the number of independent sets of size $t = \lceil 2 n\log^2\!n/d\rceil$ in the induced subgraph $G[U]$.
Then by Theorem \ref{alonrodl} and the choice of $p$,
\begin{equation}
E(|U| - |Z|) \geq pn - p^t \Bigl(\frac{2e^2 \lambda}{\log^2\!n}\Bigr)^t = pn - 1. \notag
\end{equation}
Therefore there is a set $U \subset V(G)$ such that if we remove one vertex from every independent set in $U$,
the remaining set $T$ has $|T| \geq pn - 1$ and $G[T]$ has no independent set of size $t$.
It follows that
\begin{equation}\label{dlambdabound}
r(F,t) \geq pn > \frac{n}{20\lambda} \log^2\!n. \notag
\end{equation}
This completes the proof. \qed

\section{Proof of Theorem \ref{main2}}

A key ingredient in the proof of Theorem \ref{main2} is the existence of dense bipartite graphs $G$ of high girth.
For the first statement in the theorem, we let $G$ be a bipartite graph of  girth at least twelve with parts $U$ and $V$
of sizes  $m = (q + 1)(q^8 + q^4 + 1)$ and $n = (q^3 + 1)(q^8 + q^4 + 1)$ such that every vertex of $V$ has degree $q + 1$ and every vertex
of $U$ has degree $q^3 + 1$ -- these are the incidence graphs of generalized hexagons of order $(q,q^3)$ (see~\cite{GV, vm} or~\cite{Bd}
	page 115 Corollary 5.38 for details about these constructions).

\medskip

For each $u \in U$, let $(A_u,B_u)$ be a random partition of $N_G(u)$, independently for $u \in U$. Let $H$ be the random graph with $V(H) = V$  obtained by placing a complete bipartite graph with parts $A_u$ and $B_u$
inside $N_G(u)$ for each $u \in U$. It is evident that $H$ is $C_5$-free since $G$ has girth twelve.

\medskip

Now we show every independent set in $H$ has size at most $(1 + o(1))q^8$. This is sufficient
to show $r(C_5,t) > n = (1 + o(1))t^{11/8}$. Let $I$ be a set of $t$ vertices in $H$. If $|I \cap N_G(u)| = t_u$ for $u \in U$, then
\begin{equation}
P(e(I \cap N_G(u)) = 0) = 2^{1-t_u}.
\end{equation}
Since the partitions $(A_u,B_u)$ are independent over different $u \in U$ and the sum of $t_u$ is $(q + 1)t$,
\begin{equation}
P(e(I) = 0) = \prod_{u \in U} P(e(I \cap N_G(u)) = 0) = \prod_{u \in U} 2^{1-t_u} = 2^{m - (q + 1)t}.
\end{equation}
There are ${n \choose t}$ choices of $I$, so the expected number of
independent sets of size $t$ in $H$ is
\begin{equation}
{n \choose t} 2^{m - (q + 1)t} \leq 2^{t\log_2 n + m - (q + 1)t} = 2^{m - (q + o(q))t}.
\end{equation}
Since $m = (1 + o(1))q^9$, we make take $t = (1 + o(1))q^8$ so that the above expression
decays to zero. Consequently, with high probability, every independent set in $H$ has less than $t$.

\medskip

For the second statement of Theorem \ref{main2}, the Ree-Tits octagons~\cite{GV,vm} supply requisite bipartite graphs
of girth at least sixteen -- these graphs have parts of sizes $m = (q + 1)(q^9 + q^6 + q^3 + 1)$ and $n = (q^2 + 1)(q^9 + q^6 + q^3 + 1)$ with all vertices in the larger part of
degree $q + 1$. We omit the details for this case, which are almost identical to the above. \qed

\section{Random block constructions}

Theorem \ref{main2} may be generalized as follows. Let $F$ be a graph and let $\mathcal{P} = (P_1,P_2,\dots,P_k)$   be a partition of $E(F)$ into bipartite graphs with at least one edge each. Let
$X = \{x_1,x_2,\dots,x_k\}$ be new vertices, and let $F_{\mathcal{P}}$ be the graph with $V(F_{\mathcal{P}}) = V(F) \cup X$ and edge set
\begin{equation}
E(F_{\mathcal{P}}) = \bigcup_{i = 1}^k \{\{x_i,y\} : y \in V(P_i)\}.
\end{equation}
Let $L(F)$ be the family of all graphs $F_{\mathcal{P}}$ taken over partitions $\mathcal{P}$ of $E(F)$ into paths with at least one edge each.
For instance, when $F$ is a triangle, then $L(F)$ consists of $C_4$ plus a pendant edge and $C_6$. If $F$ is a pentagon then
every member of $L(F)$ is a cycle of length at most ten plus a set of pendant edges.
Let $G$ be a bipartite graph with parts $U$ and $V$ containing no member of $L(F)$ and such that every vertex of $V$ has degree $d$.
We form a new graph $H$ with $V(H) = V$ by taking for each $u \in U$ independently a random partition $(A_u,B_u)$ of $N_G(u)$ and then adding a complete
bipartite graph with parts $A_u$ and $B_u$. By definition, $H$ does not contain $F$. Then the proof of the following is the same as the proof of Theorem \ref{main2}:

\begin{thm} \label{block}
Let $F$ be a graph and let $G$ be an $L(F)$-free bipartite graph with parts $U$ and $V$ such that $|U| = m$ and $|V| = n$ and every
vertex of $V$ has degree at least $d$. If $dt > m + t\log n$, then
\begin{equation}
r(F,t) > n.
\end{equation}
\end{thm}

If $F = K_4$, then a $C_4$-free graph containing no 1-subdivision of $K_4$ is $L(F)$-free. It is possible to show
that any graph not containing a 1-subdivision of $K_4$ has $O(n^{7/5})$ edges (see Conlon and Lee~\cite{ConlonLee}, and Janzer~\cite{Janzer}). If there is a $d$-regular
graph containing no 1-subdivision of $K_4$ with $n$ vertices and with $d = \Omega(n^{2/5})$ even, then one can produce
a random graph $H$ as above that is $d^2$-regular, and has a chance to be an $(n,d,\lambda)$-graph with $\lambda = d^{1/2 + o(1)}$
as in the work of Conlon~\cite{Conlon}. Via Theorem \ref{main}, this would then show $r(4,t) = t^{3 - o(1)}$. However, the best construction of
an $n$-vertex graph with no subdivision of $K_4$ has only $O(n^{4/3})$ edges.

\section{Concluding remarks}

$\bullet$ Although the construction in Theorem \ref{main} starting with Alon's~\cite{A}  pseudorandom triangle-free graph  provides slightly worse bounds than the known random constructions for $r(3,t)$, the number of random bits used is less than the
known random constructions~\cite{BK, S}, which use roughly $t^{4 + o(1)}$ bits.
The same observation applies to the case $r(C_{\ell}, t)$ for $\ell\in \{4,6,10\}$ where we match or exceed the best known construction obtained via the $C_{\ell}$-free process using fewer random bits.

  \medskip

  $\bullet$ It  would be interesting to see if the choice of the random subset $U$ in the proof of Theorem \ref{main} can be made explicit; for instance, the best explicit construction~\cite{KPR} of a $K_4$-free graph without independent sets of size $t$ only gives $r(4,t) = \Omega(t^{8/5})$, as compared to random graphs which give $r(4,t) = \Omega^*(t^{5/2})$.

  \medskip

  $\bullet$ If we apply the proof of Theorem \ref{main} to Paley graphs of order $q$, which are $(n,d,\lambda)$-graphs with $d = (q - 1)/2$ and $\lambda = \frac{1}{2}(\sqrt{q} \pm 1)$ where $q$ is a prime power congruent to 1 mod 4, we find almost all subsets of $\Omega(\sqrt{q}\log^2 q)$ vertices have no independent set or clique of size more than $2(\log q)^2$. In fact, Noga Alon (personal communication) had already observed a stronger statement in 1991,  that one can randomly take $q^{\alpha}$ vertices for suitable $\alpha$ and the resulting induced subgraph has clique and independence number $O(\log q)$.  It would be interesting to know if this
can be done without randomness. It is a major open question (see Croot and Lev~\cite{CL})  to determine, when $q$ is prime, the maximum size of independent sets and cliques in the Paley graph. These were shown to be at least  $\Omega(\log q \log\log q)$ by Montgomery~\cite{Mont} under GRH and at least $\Omega(\log q \log\log \log q)$ unconditionally by Graham and Ringrose~\cite{GR}. The current best upper bound is
$\sqrt{q/2} +1$ by Hanson and Petridis~\cite{HP}.
  \medskip

$\bullet$  In order to improve the exponent in the lower bound (\ref{rsn}) using Theorem \ref{main}, one could try to
find a $K_s$-free $(n,d,\lambda)$-graph with $n/\lambda \geq (n/d)^{(q + 1)/2}$ for some $q > s$, so as to obtain
$r(s,t) = \Omega(t^{(q + 1)/2})$. In the case $\lambda = O(\sqrt{d})$, it is sufficient that $d = \Omega(n^{1 - 1/q})$.

\bigskip

{\bf Acknowledgments.}
We would like to thank Noga Alon and David Conlon for helpful comments.

\end{document}